\documentclass[a4paper,fleqn]{cas-dc}

\usepackage[numbers]{natbib}
\usepackage{multirow}
\usepackage{amsmath,amssymb,amsfonts}
\usepackage{mathtools}
\usepackage{mathrsfs}
\usepackage{xcolor}
\usepackage{textcomp}
\usepackage{manyfoot}
\usepackage{booktabs}
\usepackage{algorithm}
\usepackage{algorithmicx}
\usepackage{algpseudocode}
\usepackage{listings}
\usepackage{physics}

\newcommand{\e}[1]{\begin{equation}#1\end{equation}}
\newcommand{\al}[1]{\begin{align}#1\end{align}}
\newcommand{\ald}[1]{\begin{aligned}#1\end{aligned}}
\newcommand{\ite}[1]{\begin{itemize}#1\end{itemize}}

\newcommand{\ang}[1]{\left\langle#1\right\rangle}

\def\coloneq{\mathrel{\mathop:}=}
\def\r{\mathbb{R}}
\def\rn{\mathbb{R}^n}
\def\rm{\mathbb{R}^m}
\def\sn{\mathbb{S}^n}

\newtheorem{thm}{Theorem}[section]

\newdefinition{exm}{Example}
\newdefinition{rmk}{Remark}
\newdefinition{dfn}{Definition}
\newdefinition{asm}{Assumption}
\newproof{proof}{Proof}

\begin{document}
\let\WriteBookmarks\relax
\def\floatpagepagefraction{1}
\def\textpagefraction{.001}

\shorttitle{Nonconvex SDP with globally optimal KKT points}
\shortauthors{A. Nishioka and Y. Kanno}

\title[mode = title]{A class of nonconvex semidefinite programming in which every KKT point is globally optimal}

\author[1]{Akatsuki Nishioka}
\cormark[1]
\ead{nishioka.a.2122@m.isct.ac.jp}
\credit{Conceptualization of this study, Methodology, Writing}

\author[2,3]{Yoshihiro Kanno}
\credit{Finding applications, Instruction}

\affiliation[1]{organization={Department of Mathematical and Computing Science, School of Computing, Institute of Science Tokyo},
                addressline={Ookayama 2-12-1},
                city={Meguro-ku},
                postcode={152-8550},
                state={Tokyo},
                country={Japan}}

\affiliation[2]{organization={Department of Mathematical Informatics, Institute of Science and Technology, The University of Tokyo},
                addressline={Hongo 7-3-1},
                city={Bunkyo-ku},
                postcode={113-8656},
                state={Tokyo},
                country={Japan}}

\affiliation[3]{organization={Mathematics and Informatics Center, The University of Tokyo},
                addressline={Hongo 7-3-1},
                city={Bunkyo-ku},
                postcode={113-8656},
                state={Tokyo},
                country={Japan}}

\cortext[cor1]{Corresponding author}

\begin{abstract}
We consider a special class of nonconvex semidefinite programming problems and show that every point satisfying the Karush--Kuhn--Tucker (KKT) conditions is globally optimal despite nonconvexity. This property is related to pseudoconvex optimization and fractional programming. We also present several applications to robust fractional programming and generalized eigenvalue optimization appearing in topology optimization, network control, finance, etc. 
\end{abstract}

\begin{keywords}
semidefinite programming \sep pseudoconvex optimization \sep fractional programming
\MSC[2020]{90C22 \sep 90C26}
\end{keywords}

\maketitle

\section{Introduction}

There is a long history of identifying nonconvex optimization problems
in which first-order stationarity guarantees global optimality.
Pseudoconvex and invex functions are prominent examples, and
fractional programming is closely related to these notions under
suitable assumptions
\cite{avriel10,mishra08,stancu12}.
By contrast, results establishing the global optimality of every KKT
point in nonconvex nonlinear semidefinite programming (NSDP) remain
scarce. Although invexity has been studied for NSDP \cite{sun11}, it is generally
difficult to verify because its definition involves an unknown kernel
function.

The present study is motivated by a generalized eigenvalue optimization arising in truss topology optimization in structural engineering. This problem was previously shown to be pseudoconvex minimization, and hence every Clarke stationary point is globally optimal
\cite{nishioka23coap}. In this paper, rather than focusing on this particular problem, we seek a broader NSDP framework under which the same global-optimality property holds.

The main contributions of this paper are as follows.
\begin{itemize}
  \item
  We introduce a new class of nonconvex SDPs and prove that every KKT
  point is globally optimal.
  \item
  The proposed framework extends the generalized eigenvalue
  optimization setting studied in \cite{nishioka23coap} and covers new
  applications, including robust fractional programming and eigenvalue optimization in network control.
  \item
  The NSDP formulation can accommodate more complicated constraints
  than problem-specific formulations. It also allows the problems to be treated in a unified manner using open-source NSDP solvers, such as PENLAB
  \cite{fiala13penlab} and IPNSDP \cite{weldeyesus26}.
  \item
  Since every KKT point is globally optimal, problems in the proposed
  class provide convenient benchmark instances for nonconvex SDP
  algorithms.
  \item We show that several classes of pseudoconvex functions and optimization problems admit convenient semidefinite-programming representations.
  \item We investigate the relationship between the proposed problem class
  and pseudoconvex optimization. In particular, we establish the
  quasiconvexity of an associated value function and discuss its
  connection to pseudoconvexity.
\end{itemize}

The remainder of this paper is organized as follows.
Section~\ref{sec_prob} introduces the problem class and its
assumptions. Section~\ref{sec_app} presents several
applications. Section~\ref{sec_glob} proves the global optimality of
KKT points. Section~\ref{sec_rel} discusses the relationship
with pseudoconvex optimization, and
Section~\ref{sec_conc} concludes the paper.

\subsection*{Notation and definitions}\label{sec_not}

We use the following notation:
\ite{
\item $\sn$ is the set of $n\times n$ symmetric matrices.
\item $A\succeq0$ means that $A$ is a symmetric positive semidefinite matrix. $A\succ0$ means that $A$ is a symmetric positive definite matrix.
\item $\ang{A,B}\coloneq\mathrm{tr}(AB)$ is the inner product of symmetric matrices $A,B\in\sn$. We often write $\ang{A,B}_{\sn}$ and $\ang{x,y}_{\rm}$ to distinguish between these two inner products.
\item For $F:\rm\to\sn$, $\nabla F(x)\in\sn\times\rm$ is the Jacobian of $F$ at $x\in\rm$, which is an array consisting of $\frac{\partial}{\partial x_i} F(x) \in\sn\ (i=1,\ldots,m)$. For $A:\rm\times\r\to\sn$, $\nabla_x A(x,y)\in\sn\times\rm$ is the Jacobian of $A$ with respect to the first variable $x\in\rm$ at $(x,y)$. $\frac{\partial}{\partial y}A(x,y)\in\sn$ is the partial derivative of $A$ with respect to the second variable $y\in\r$ at $(x,y)$.
}

We also use the following definition of convexity and concavity of a symmetric-matrix-valued function.

\begin{dfn}[matrix convexity and concavity]\label{dfn_conv_mat}
We say that a symmetric-matrix-valued function $F:\rm\to\sn$ is convex if
\e{\label{conv_mat}
\lambda F(x)+(1-\lambda)F(x')\succeq F(\lambda x+(1-\lambda)x')
}
holds for any $x,x'\in\rm$ and $\lambda\in[0,1]$. We also say $F$ is concave if $-F$ is convex.
% We say $F$ is strictly convex if \eqref{conv_mat} is satisfied with $\succ$ instead of $\succeq$ for any $x,x'\in\rm$ such that $x\neq x'$.
%  
This definition is taken from \cite[Definition 2.3.2]{wolkowicz12}.
\end{dfn}

Since $A\succeq0$ if and only if $\ang{A,Z}\ge0$ for any $Z\succeq0$ \cite[Theorem 2.3.11]{wolkowicz12}, the above definition implies that $F:\rm\to\sn$ is convex if and only if $x\mapsto\ang{F(x),Z}$ is convex for any $Z\succeq0$. If $F:\rm\to\sn$ is convex and differentiable, then we have, for any fixed $Z\succeq0$,
\e{\label{conv_mat_diff}
\ang{F(y),Z}_{\sn}-\ang{F(x),Z}_{\sn}\ge \ang{\ang{\nabla F(x),Z}_{\sn},y-x}_{\rm}
}
for any $x,y\in\rm$, where $\nabla F(x)\in\sn\times\rm$ is the Jacobian of $F(x)$. We use \eqref{conv_mat_diff} later.

\section{Problem setting}\label{sec_prob}

We consider the following class of nonconvex semidefinite programming (SDP) problems:
\e{
\ald{
& \underset{x\in\rm,\ y\in\r}{\mathrm{minimize}} & & y \\
& \mathrm{subject\ to} & & A(x,y)\succeq0,\\
& & & B(x)\succeq0,
}
\label{p}
\tag{P}
}
with the following assumptions.
\begin{asm}\label{asm1}
\ \\ \vspace{-6mm}
\begin{enumerate}[(a)]
\item $A(\cdot,y):\rm\to\sn$ is differentiable and concave for any feasible $y\in\r$. $B(\cdot):\rm\to\sn$ is differentiable and concave (see Definition \ref{dfn_conv_mat}).
\item $A(x,\cdot):\r\to\sn$ is differentiable and strictly increasing in the Loewner order $\succ$ for every feasible $x$; that is, $y>y'$ implies $A(x,y)\succ A(x,y')$ (a sufficient condition is $\frac{\partial}{\partial y}A(x,y)\succ0$).
\end{enumerate}
\end{asm}

The restriction that the objective in Problem \eqref{p} is the scalar variable
$y$ does not cause an essential loss of generality. Indeed, a nonlinear SDP
\e{
\begin{aligned}
\underset{x\in\mathbb{R}^m}{\mathrm{minimize}}
  \quad & f(x) \\
  \mathrm{subject\ to}
  \quad & B(x)\succeq 0
\end{aligned}
}
can be equivalently written in the epigraph form
\e{
\begin{aligned}
\underset{x\in\mathbb{R}^m,\ y\in\mathbb{R}}{\mathrm{minimize}}
  \quad & y \\
  \mathrm{subject\ to}
  \quad & y-f(x)\geq 0, \\
        & B(x)\succeq 0
\end{aligned}
}
by introducing the auxiliary variable $y$. Assumption \ref{asm1} does impose restrictions, but, as shown in Section \ref{sec_app}, the resulting
problem class still covers several important applications. For example, some convex SDP problems with $f$ convex and $B$ concave can be reformulated as Problem \eqref{p} satisfying Assumption \ref{asm1}.

Note that the
constraint $B(x)\succeq 0$ can incorporate scalar inequalities such as
$g(x)\geq 0$ with a concave function $g$ as diagonal blocks. The constraint $B(x)\succeq 0$ cannot be included in $A(x,y)\succeq0$ since Assumption \ref{asm1}(b) no longer holds. Also, matrix optimization variables can be accommodated by vectorizing them into $x$.

For simplicity, affine equality constraints appearing in the applications below are omitted in Problem \eqref{p}. Such additional constraints do not change the results of this paper. Note that if we treat an affine equality constraint as two affine inequality constraints and include them in $B(x)\succeq0$, it can break a constraint qualification.

\section{Applications}\label{sec_app}

The main applications of Problem \eqref{p} have the following form of nonsmooth fractional programming
\e{
\underset{x\in \mathcal{C}\subseteq\rm}{\mathrm{minimize}}\ \underset{q\in \mathcal{Q}}{\sup}\ \frac{f(x,q)}{g(x,q)},
}
where $\mathcal{C}$ is the convex feasible set, $\mathcal{Q}$ is the (finite or infinite) uncertainty set, and $g(x,q)>0$ for any $x\in \mathcal{C}$ and $q\in \mathcal{Q}$. This problem can be reformulated as
\e{
\ald{
& \underset{x\in\mathcal{C}\subseteq\rm,\ y\in\r}{\mathrm{minimize}} & & y \\
& \mathrm{subject\ to} & & yg(x,q)-f(x,q)\ge 0\ \ \forall q\in \mathcal{Q}.\\
}
}
If these (possibly semi-infinite) constraints admit a semidefinite representation and satisfy Assumption \ref{asm1}, then this problem can be written as Problem \eqref{p}. Problems of this type arise in generalized eigenvalue optimization and robust fractional programming. In the following, we introduce several applications. For other applications of robust fractional programming in operations research, see \cite{gorissen15}.

\subsection{Eigenfrequency topology optimization}\label{sec_to}

Topology optimization is a method to find an optimal design of a structure and has a wide variety of applications in engineering \cite{bendsoe04}.
We consider a maximization problem of the smallest generalized eigenvalue (the squared fundamental frequency) of a truss structure in topology optimization. By maximizing this quantity, we can obtain a structure resistant to vibration. See \cite{achtziger07siam,achtziger07smo,nishioka23coap,nishioka24svva,ohsaki99} for more details on eigenfrequency topology optimization.

The problem is formulated as a maximization problem of the minimum generalized eigenvalue
\e{
\ald{
& \underset{x\in[x_{\min},\infty)^m}{\mathrm{maximize}} & & \lambda_\mathrm{min}(K(x),M(x))\\
& \mathrm{subject\ to} & & l^\top x \le V_0,\\
}
\label{p_to}
}
where $x\in[x_{\min},\infty)^m$ is a vector consisting of cross-sectional areas of bars of the truss structure, $m$ is the number of bars, $x_{\min}>0$ is the lower bound of the cross-sectional areas\footnote{For studies of cases with $x_\mathrm{min}=0$, see \cite{achtziger07siam,nishioka24svva}.}, $K(x)=\sum_{j=1}^m x_j K_j$ and $M(x)=M_0+\sum_{j=1}^m x_j M_j$, where $K_j,M_j\in\sn$ $(j=1,\ldots,m)$ are the element stiffness and mass matrices, respectively, $M_0\in\sn$ is a nonstructural mass matrix, and $K_j\succeq0$, $M_j\succeq0$, and $M_0\succeq0$. We assume that $K(x)$ and $M(x)$ are symmetric positive definite for any $x\in(0,\infty)^m$ after imposing appropriate boundary conditions. Moreover, $l\in(0,\infty)^m$ is a vector consisting of undeformed lengths of bars of a truss structure, and $V_0>0$ is the upper bound of the volume of the truss structure. The minimum generalized eigenvalue is defined by
\e{\ald{
&\lambda_\mathrm{min}(K(x),M(x))\\
& \coloneq \min\{\lambda_i\in\r\ (i=1,\ldots,n)\mid\\
&\qquad\qquad\exists v_i\in\rn\backslash\{0\}\text{ s.t. }K(x)v_i=\lambda_i M(x)v_i\}\\
& = \underset{v\in\rn\backslash\{0\}}{\inf}\frac{v^\top K(x)v}{v^\top M(x)v}\\
& = \sup\{\lambda\in\r\mid K(x)-\lambda M(x)\succeq0\}.
\label{geneig}
}}
Problem \eqref{p_to} is shown to be a pseudoconcave maximization problem in \cite{nishioka23coap}, and thus its stationary points (in the sense of Clarke) are always globally optimal.

By using the last formulation of the minimum generalized eigenvalue in \eqref{geneig}, Problem \eqref{p_to} can be reformulated into a nonconvex SDP
\e{
\ald{
& \underset{x\in[x_{\min},\infty)^m,\ \lambda\in\r}{\mathrm{minimize}} & & -\lambda\\
& \ \ \ \mathrm{subject\ to} & & K(x)-\lambda M(x)\succeq 0,\\
& & & l^\top x\le V_0.
}
\label{p_to_nsdp}
}
By setting $y=-\lambda$, Problem \eqref{p_to_nsdp} satisfies Assumption~\ref{asm1}. In \cite{nishioka23coap}, a smoothing gradient method was proposed for Problem~\eqref{p_to}. However, if more complicated constraints are added to Problem~\eqref{p_to}, the smoothing gradient method cannot be applied. On the other hand, the NSDP approach is more convenient for treating complicated constraints.

\subsection{Grasping force optimization in robotics}\label{sec_robo}

We consider the following nonlinear second-order cone programming problem
\e{
\ald{
& \underset{x_i\in\rm,\ y\in\r}{\mathrm{minimize}} & & y \\
& \mathrm{subject\ to} & & \|A_i x_i\|\le y b_i^\top x_i\ \ (i=1,\ldots,n),\\
& & & b_i^\top x_i \ge \beta>0\ \ (i=1,\ldots,n),
}
\label{p_nsocp}
}
where $A_i\in\r^{l\times m}$ and $b_i\in\rm$ are constant matrices and vectors. The nonlinearity appears in the multiplication $y b_i^\top x_i$ of two variables $x_i$ and $y$. The nonlinear second-order cone constraints can be written by nonlinear semidefinite constraints
\e{
\begin{pmatrix}
y b_i^\top x_i I & A_i x_i \\
x_i^\top A_i^\top & y b_i^\top x_i
\end{pmatrix}
\succeq 0\ (i=1,\ldots,n),
}
which, together with $b_i^\top x_i\ge\beta>0$, satisfy Assumption~\ref{asm1}.

A possible application is a modified version of grasping force optimization in robotics \cite{lobo98}. We consider a situation where a robot hand tries to grasp an object (rigid body with smooth boundary) with $n$ fingers. In this setting, $x_i\in\r^3$ is the force of the $i$-th finger, $b_i$ is the unit normal vector to the surface of the object at the contact point of the $i$-th finger, $y$ is the friction coefficient, and $A_i=I- b_i b_i^\top$ is the projection matrix onto the tangent plane of the contact surface. Additionally, we consider equilibrium equations:
\e{
\sum_{i=1}^n x_i + f_\mathrm{ext} = 0,\ \ \sum_{i=1}^n p_i \times x_i + T_\mathrm{ext} = 0
\label{equil}
}
where $f_\mathrm{ext}\in\r^3$ is the external force applied to the object, $p_i\in\r^3$ is the position of the $i$-th contact point, $\times$ is the vector product, and $T_\mathrm{ext}\in\r^3$ is the external torque. Problem \eqref{p_nsocp} with constraints \eqref{equil} seeks to find the minimum friction coefficient for which the robot hand can grasp the object with suitable choices of forces $x_i$.

Note that Problem \eqref{p_nsocp} can also be written as the following nonsmooth fractional programming problem:
\e{
\underset{x_i\in \rm}{\mathrm{minimize}}\ \underset{i}{\max}\frac{\|A_i x_i\|}{b_i^\top x_i}
}
subject to $b_i^\top x_i\ge\beta$ $(i=1,\ldots,n)$ and the constraints \eqref{equil}.

\subsection{Robust portfolio optimization}

Portfolio optimization is a basic problem in finance. We consider a robust counterpart of a portfolio optimization problem (inverse Sharpe ratio minimization) \cite{palomar25}
\e{\ald{\label{eq:worst_case_ir}
& \underset{x\in\rm}{\mathrm{minimize}} & & \max_{i=1,\ldots,n} \frac{\sqrt{x^\top \Sigma_i x}}
{\mu_i^\top x} \\
& \mathrm{subject\ to} & &  \mu_i^\top x\ge\beta>0\ \ (i=1,\ldots,n)\\
& & & e^\top x=1,\ 0\le x\le u
}}
where $x\in\rm$ is the portfolio vector,
$\Sigma_i \succeq 0$ is the covariance matrix in regime $i$, $\mu_i\in\rm$ is the
expected excess-return vector in regime $i$, $e\in\rm$ is the all-ones vector, $u\in\rm$ is the upper bound of $x$. Since $\sqrt{x^\top \Sigma_i x}=\|\Sigma_i^{1/2}x\|$, using the same argument as in Section \ref{sec_robo}, we can reformulate this problem into the form of Problem \eqref{p} satisfying Assumption~\ref{asm1}.

\subsection{Optimization of connectivity in network control}

We consider a connected weighted graph
$\mathcal{G}=(\mathcal{N},\mathcal{E})$ representing a power
transmission network, as studied by Guo et al.~\cite{guo18}.
The nodes and edges represent buses and transmission lines,
respectively. Let $C\in\mathbb{R}^{n\times m}$ be an incidence matrix.
More precisely, the column $c_e$ of $C$
corresponding to an edge $e=(i,j)$ has entry $1$ at its source
node $i$, entry $-1$ at its terminal node $j$, and zero entries elsewhere. Moreover, let
\e{
M=\operatorname{diag}(M_1,\ldots,M_n)\succ0,\quad
B=\operatorname{diag}(B_e:e\in\mathcal{E})\succ0
}
denote the inertia matrix and the line-susceptance matrix, respectively. The weighted graph Laplacian is $CBC^\top$, and
Guo et al.~\cite{guo18} consider its inertia-scaled form
\e{
L\coloneqq M^{-1/2}CBC^\top M^{-1/2}.
}
This matrix is positive semidefinite and $M^{1/2}\mathbf{1}$ is an eigenvector corresponding to the smallest eigenvalue $\lambda_1(L)=0$.

For a connected network, its second smallest eigenvalue
$\lambda_2(L)$ measures
connectivity of the scaled network; larger $\lambda_2(L)$ suggests stronger connectivity \cite{guo18}. Let $V\in\mathbb{R}^{n\times(n-1)}$ be a full-column-rank matrix whose columns span the subspace $\{q\in\mathbb{R}^n:q^\top M\mathbf{1}=0\}$. The Rayleigh quotient characterization gives
\begin{align}
    \lambda_2(L)
    &=\min_{\substack{q\neq0\\q^\top M\mathbf{1}=0}}
    \frac{q^\top CBC^\top q}{q^\top Mq} \notag\\
    &=
    \lambda_{\min}
    \left(
    V^\top CBC^\top V,\,
    V^\top MV
    \right).
    \label{eq:scaled-algebraic-connectivity}
\end{align}

Although Guo et al.~\cite{guo18} only consider analysis of $\lambda_2(L)$, we can use their model to formulate a network-design problem. To keep $V$ independent of the design variables, let $M=\alpha M_0$, where $M_0\succ0$ is fixed, $\alpha>0$, and $V^\top M_0\mathbf{1}=0$. Let $b\in\mathbb{R}^m$ collect controllable line susceptances and set $B(b)=\operatorname{diag}(b)$. Then the maximization of $\lambda_2(L)$ can be written as
\e{\ald{
    & \underset{b\in\mathbb{R}^m,\ \alpha\in\mathbb{R},\ y\in\r}{\operatorname{minimize}}
    \quad& &
    y\\
    & \text{subject to}\quad& &
    V^\top C B(b) C^\top V
    +
    y\alpha V^\top M_0 V
    \succeq0,\\
    &&&(b,\alpha)\in\mathcal{C},
}}
where $\mathcal{C}$ is a compact convex set encoding, for example, bounds and budget constraints, with $B(b)\succeq\underline b I$ for some $\underline b>0$ and $\alpha\ge\underline\alpha>0$. For every fixed $y$, the matrix inequality is affine in $(b,\alpha)$, and its derivative with respect to $y$ is $\alpha V^\top M_0V\succ0$. Thus, the problem satisfies Assumption~\ref{asm1}, while the product $y\alpha$ makes it generally nonconvex. We emphasize that the assumption $M=\alpha M_0$ is rather restrictive. For a general matrix variable $M$, the resulting formulation is a more general NSDP and may no longer satisfy Assumption~\ref{asm1}. Determining when network-design problems still enjoy global optimality of KKT points is left for future work.

The second smallest eigenvalue of the graph Laplacian is often called the algebraic connectivity. For other forms of optimization of the algebraic connectivity, see also \cite{ghosh06,kim06}.

\section{Global optimality of KKT points}\label{sec_glob}

Under a certain constraint qualification (e.g.~the Mangasarian--Fromovitz constraint qualification), the KKT condition provides a necessary condition for local optimality \cite{yamashita15}. The KKT condition for Problem \eqref{p} is as follows:
\begin{subequations}\label{kkt}
\al{
& \ang{\nabla_x A(x,y),Z}+\ang{\nabla B(x),W}=0\label{kkt1}\\
& 1-\ang{\pdv{y}A(x,y),Z}=0,\label{kkt2}\\
& \ang{A(x,y),Z}=0,\ \ang{B(x),W}=0,\label{kkt3}\\
& A(x,y)\succeq0,\ B(x)\succeq0,\ Z\succeq0,\ W\succeq0,\label{kkt4}
}
\end{subequations}
where $Z,W\in\sn$ are the Lagrange multipliers and $\nabla_x A(x,y)\in\sn\times\rm$ is the Jacobian of $A(x,y)$ with respect to $x$ (note that $\ang{\nabla_x A(x,y),Z}\in\rm$ and the first equality is taken in $\rm$). 

Although the KKT condition is only a necessary condition for global optimality in general nonconvex SDP, we can show that it is also sufficient in Problem \eqref{p}.

\begin{thm}\label{thm}
Under Assumption~\ref{asm1}, the primal part $(x^*,y^*)$ of every KKT tuple $(x^*,y^*,Z^*,W^*)$ satisfying \eqref{kkt} is globally optimal for Problem \eqref{p}.
\end{thm}
\begin{proof}
Let $(x^*,y^*,Z^*,W^*)$ be a KKT tuple, and let $(\bar{x},\bar{y})$ be an arbitrary feasible point. Suppose, to the contrary, that $\bar{y}<y^*$.

Since $A(\cdot,y^*)$ and $B$ are concave and $Z^*,W^*\succeq0$, the scalar-valued functions $x\mapsto\ang{A(x,y^*),Z^*}_{\sn}$ and $x\mapsto\ang{B(x),W^*}_{\sn}$ are concave. Hence,
\e{\ald{\label{eq:proof}
&\ang{A(\bar{x},y^*),Z^*}_{\sn}+\ang{B(\bar{x}),W^*}_{\sn}\\
&\le
\ang{A(x^*,y^*),Z^*}_{\sn}+\ang{B(x^*),W^*}_{\sn}\\
&\quad+
\ang{
\ang{\nabla_x A(x^*,y^*),Z^*}_{\sn}
+
\ang{\nabla B(x^*),W^*}_{\sn},
\bar{x}-x^*
}_{\rm}\\
&=0,
}}
where the equality follows from stationarity~\eqref{kkt1} and complementarity~\eqref{kkt3}.

On the other hand, feasibility of $(\bar{x},\bar{y})$ and Assumption~\ref{asm1}(b) give
\e{
A(\bar{x},y^*)\succ A(\bar{x},\bar{y})\succeq0.
}
Also, $Z^*\neq0$ because~\eqref{kkt2}:
\e{
\ang{\pdv{y}A(x^*,y^*),Z^*}=1.
}
Therefore,
\e{
\ang{A(\bar{x},y^*),Z^*}>0.
}
Moreover, since $B(\bar{x})\succeq0$ and $W^*\succeq0$, we have $\ang{B(\bar{x}),W^*}\ge0$. These contradict the inequality \eqref{eq:proof}. Thus, no feasible point has an objective value smaller than $y^*$, and $(x^*,y^*)$ is globally optimal.\qed
\end{proof}

\begin{rmk}
Neither a constraint qualification nor the existence of a global
minimizer is required in Theorem~\ref{thm}. A constraint
qualification is needed only to ensure that a local minimizer satisfies
the KKT conditions; in that case, every local minimizer is globally
optimal by Theorem~\ref{thm}. An existence assumption may be
imposed separately when one wishes to guarantee that
Problem~\eqref{p} admits an optimal solution.
\end{rmk}

% With the strict concavity assumption on $A$, we obtain the uniqueness of solutions.

% \begin{prp}
% In addition to Assumption \ref{asm1}, we assume that $A(\cdot,y):\rm\to\sn$ is strictly concave. Then, there exists the unique KKT point, which is globally optimal.
% \end{prp}
% \begin{proof}
% Suppose there exist two KKT tuples $(x^*,y^*,Z^*,W^*)$ and $(x',y',Z',W')$ such that $x^*\neq x'$ and $y^*=y'$. Then, by the same argument as \eqref{eq3} with strict inequality, we obtain
% \e{
% 0>\ang{A(x',y^*),Z^*}+\ang{B(x'),W^*}.
% }
% In contrast, by the same argument as \eqref{eq4}, we obtain 
% \e{
% \ang{A(x',y^*),Z^*}\ge\ang{A(x',y'),Z^*}+\ang{\pdv{y}A(x',y'),Z^*}(y^*-y')=\ang{A(x',y'),Z^*}\ge 0.
% }
% We also have $\ang{B(x'),W^*}\ge0$ by \eqref{kkt4}. Thus, it results in a contradiction.
% \end{proof}

\section{Relation to pseudoconvexity}\label{sec_rel}

We first introduce the definitions of pseudoconvex functions and quasiconvex functions. It is known that every pseudoconvex function is quasiconvex \cite{penot97}.

\begin{dfn}[pseudoconvex function \cite{penot97}]
A differentiable function $f:\rm\to\r$ is said to be pseudoconvex if the following holds:
\e{
\forall x,y\in\rm,\ f(x)>f(y)\ \Rightarrow\ \ang{\nabla f(x),y-x}<0.
\label{pseudo}
}
\end{dfn}

For the definition of pseudoconvex functions for nonsmooth functions, see \cite{penot97}.

\begin{dfn}[quasiconvex function \cite{penot97}]
A function $f:\rm\to(-\infty,\infty]$ is said to be quasiconvex if its sublevel set
\e{
\{x\in\rm \mid f(x)\le\alpha\}
}
is convex for any $\alpha\in\mathbb{R}$.
\end{dfn}

The problem \eqref{p} satisfying Assumption \ref{asm1} is equivalently written as an unconstrained minimization problem of the function $\varphi:\rm\to(-\infty,\infty]$ defined by
\e{\label{value}
\varphi(x)\coloneqq\inf\{y\in\r\mid A(x,y)\succeq 0,\ B(x)\succeq 0\}.
}
Here, $\inf\emptyset\coloneqq+\infty$.
The following result shows $\varphi$ is quasiconvex. 

\begin{thm}
Under Assumption \ref{asm1}, the function $\varphi$ defined by \eqref{value} is quasiconvex.
\end{thm}
\begin{proof}
We show that, for any $\alpha\in\r$, the sublevel set of $\varphi$ is
\e{
\{x\in\rm\mid\varphi(x)\le\alpha\}
=\{x\in\rm\mid A(x,\alpha)\succeq0,\ B(x)\succeq0\}.
}

Suppose $A(x,\alpha)\succeq0,\ B(x)\succeq0$. Then $y=\alpha$ is feasible in the definition of $\varphi(x)$, and hence $\varphi(x)\leq\alpha$. This proves the inclusion ``$\supseteq$''.

Conversely, suppose that $\varphi(x)\leq\alpha$.
If $\varphi(x)<\alpha$, the definition of the infimum implies
that there exists a feasible $y<\alpha$ such that $A(x,y)\succeq0,\ B(x)\succeq0$. By Assumption~\ref{asm1}(b),
\e{
A(x,\alpha)\succ A(x,y)\succeq0.
}

It remains to consider the case $\varphi(x)=\alpha$.
By the definition of the infimum, there exists a sequence
$\{y_k\}$ such that
\e{
A(x,y_k)\succeq0,\quad B(x)\succeq0,\quad y_k\to\alpha.
}
Since $A(x,\cdot)$ is continuous and the positive semidefinite
cone is closed, we obtain
\[
    A(x,\alpha)
    =
    \lim_{k\to\infty}A(x,y_k)
    \succeq0.
\]
This proves the inclusion ``$\subseteq$''.

Since
$A(\cdot,\alpha)$ and $B$ are concave, this sublevel set is convex.
Hence, $\varphi$ is quasiconvex.\qed
\end{proof}

Analyzing pseudoconvexity of $\varphi$ is much more difficult because it involves the subdifferential $\partial\varphi$, i.e., sensitivity analysis of nonconvex SDP $\min\{y\in\r\mid A(x,y)\succeq 0,\ B(x)\succeq 0\}$ with respect to $x$ (cf.~\cite[Section 4.1.4]{wolkowicz12}).

However, we can see that $\varphi$ behaves like a pseudoconvex function by the following formal calculation. To simplify the discussion, let us ignore the constraint $B(x)\succeq0$. Suppose $\varphi$ is differentiable at $x$ and the regularity conditions required for the sensitivity formula and uniqueness of the Lagrange multiplier $Z(x)\succeq0$ hold (see \cite[Chapter 4]{wolkowicz12} for details). The Lagrange multiplier satisfies
\e{\ald{
&A(x,\varphi(x))\succeq 0,\ 
\langle A(x,\varphi(x)), Z(x)\rangle = 0,\\
&1-
\ang{\frac{\partial}{\partial y} A(x,\varphi(x)),Z(x)}=0.
}}
Then, \cite[Theorem 4.1.11]{wolkowicz12} yields
\e{
\nabla\varphi(x)=-\langle \nabla_x A(x,\varphi(x)), Z(x)\rangle.
}

Consider an open set $U\subseteq \mathrm{dom}\,\varphi$ where $\varphi$ is differentiable. Take $x,x'\in U$ and suppose that $\varphi(x')<\varphi(x)$. Since $A(x',\varphi(x'))\succeq 0$, Assumption \ref{asm1}(b) implies $A(x',\varphi(x))\succ 0$. On the other hand, the concavity of $x\mapsto A(x,y)$ gives
\e{\ald{
&\langle A(x',\varphi(x)), Z(x)\rangle\\
&\le
\langle A(x,\varphi(x)), Z(x)\rangle
+
\ang{\langle \nabla_x A(x,\varphi(x)), Z(x)\rangle,x'-x}.
\label{pr_1}
}}
By complementarity, the first term on the right-hand side is $0$. Moreover, $Z(x)\neq 0$ because
$\langle \partial_y A(x,\varphi(x)), Z(x)\rangle = 1$, and therefore $A(x',\varphi(x))\succ0$ implies $\langle A(x',\varphi(x)), Z(x)\rangle>0$. Above all, \eqref{pr_1} implies 
\e{
\ang{\langle \nabla_x A(x,\varphi(x)), Z(x)\rangle,x'-x}>0,
}
which is equivalent to
\e{
\ang{\nabla \varphi(x),x'-x}<0.
}
This is exactly the pseudoconvex implication.

We emphasize that the argument above is only a proof strategy under strong assumptions. Extending it to the general nonsmooth case requires a much more careful subdifferential and sensitivity analysis. This is left as a topic for future research.

\section{Conclusion}\label{sec_conc}

We proved the global optimality of KKT points in a certain class of nonconvex SDPs. Such a result appears to be new in the NSDP literature.

Assumption \ref{asm1} is restrictive, although it includes several important applications. We hope to broaden the class of nonconvex SDPs for which every KKT point is globally optimal. The relationship between Problem~\eqref{p} and pseudoconvex
(or invex) optimization should also be studied more thoroughly. It would be useful if we could develop a characterization of invex and pseudoconvex optimization problems in terms of nonconvex SDP formulations.

\section*{Acknowledgements}

The work of the first author is partially supported by JSPS KAKENHI JP25KJ0120. The work of the second author is partially supported by JSPS KAKENHI JP21K04351 and JP24K07747.

\section*{Data availability}

This study did not use any datasets.

\section*{Declaration of competing interest}

The authors declare that they have no conflict of interest.

\section*{Declaration of generative AI and AI-assisted technologies in the writing process}

During the preparation of this work, the authors used ChatGPT (GPT-5.5) to improve the clarity of the exposition, refine the presentation of mathematical proofs, and explore possible applications. After using this tool, the authors reviewed, verified, and edited the content as needed and take full responsibility for the content of the publication.

\printcredits

\bibliographystyle{cas-model2-names}
\bibliography{ref}

\end{document}